\pdfoutput=1
\RequirePackage{ifpdf}
\ifpdf 
\documentclass[pdftex]{sigma}
\else
\documentclass{sigma}
\fi

\numberwithin{equation}{section}
\newtheorem{Theorem}{Theorem}[section]
\newtheorem{Lemma}[Theorem]{Lemma}
\newtheorem{Proposition}[Theorem]{Proposition}
{ \theoremstyle{definition}
\newtheorem{Definition}[Theorem]{Definition}
\newtheorem{Remark}[Theorem]{Remark} }

\DeclareMathOperator{\spn}{span}
\DeclareMathOperator{\Ran}{Ran}

\begin{document}

\newcommand{\arXivNumber}{1410.0733}

\allowdisplaybreaks

\renewcommand{\PaperNumber}{012}

\FirstPageHeading

\ShortArticleName{The Quantum Pair of Pants}

\ArticleName{The Quantum Pair of Pants}

\Author{Slawomir KLIMEK~$^\dag$, Matt MCBRIDE~$^\ddag$, Sumedha RATHNAYAKE~$^\dag$ and Kaoru SAKAI~$^\dag$}

\AuthorNameForHeading{S.~Klimek, M.~McBride, S.~Rathnayake and K.~Sakai}

\Address{$^\dag$~Department of Mathematical Sciences, Indiana University-Purdue University Indianapolis,\\
\hphantom{$^\dag$}~402 N.~Blackford St., Indianapolis, IN 46202, USA}
\EmailD{\href{mailto:sklimek@math.iupui.edu}{sklimek@math.iupui.edu},
\href{mailto:srathnay@iupui.edu}{srathnay@iupui.edu}, \href{mailto:ksakai@iupui.edu}{ksakai@iupui.edu}}

\Address{$^\ddag$~Department of Mathematics, University of Oklahoma, 601 Elm St., Norman, OK 73019, USA}
\EmailD{\href{mailto:mmcbride@math.ou.edu}{mmcbride@math.ou.edu}}

\ArticleDates{Received October 24, 2014, in f\/inal form February 03, 2015; Published online February 10, 2015}

\Abstract{We compute the spectrum of the operator of multiplication by the complex coordinate in a~Hilbert space of
holomorphic functions on a~disk with two circular holes.
Additionally we determine the structure of the $C^*$-algebra generated by that operator.
The algebra can be considered as the quantum pair of pants.}


\Classification{46L35}

\section{Introduction}
In this paper we study the operator~$z$ of multiplication by the complex coordinate in Hilbert spaces of holomorphic
functions on certain multiply connected domains in the complex plane.
The domains we consider are disks with circular holes.
The case of a~disk with no holes is the classical one.
In the Hardy space of the disk the multiplication operator~$z$ is the unilateral shift whose spectrum is the disk.
The $C^*$-algebra generated by the unilateral shift, the Toeplitz algebra, is an extension of the algebra of compact
operators by $C(S^1)$, $S^1$ being the boundary of the disk~\cite{Cob}.
For the Bergman space the operator~$z$ is a~weighted unilateral shift and its spectrum and the $C^*$-algebra it
generates are the same as in the Hardy space~\cite{KL}.
Partially for those reasons the Toeplitz algebra is often considered as the quantum disk~\cite{KL,KL2, KL3}.

A disk with one hole is biholomorphic to an annulus.
In the Bergman space for example, the~$z$ operator is a~weighted bilateral shift with respect to the natural basis of
(normalized) powers of the complex coordinate.
Its spectrum is the annulus, and the $C^*$-algebra it generates is an extension of the algebra of compact operators~by
$C(S^1\times S^1)$, where $S^1\times S^1$ is the boundary of the annulus.
The same is true for many other Hilbert spaces of holomorphic functions on an annulus.
The resulting $C^*$-algebra is the quantum annulus of~\cite{KL3,KM1}.

In this paper we study in detail the two hole case: a~pair of pants.
Up to biholomorphism we can realize a~disk with two holes as an annulus centered at zero with outer radius one, with an
additional of\/f centered hole.
In the space of continuous functions on the closed pair of pants that are holomorphic in its interior, we consider
a~specif\/ic inner product with respect to which the operator of multiplication by the complex coordinate~$z$ has
a~particularly simple structure.
The results we obtain are completely analogous to zero and one-hole cases: the spectrum of~$z$ is the domain of the
corresponding pair of pants while the $C^*$-algebra generated by~$z$ is an extension of the algebra of compact operators
by $C(S^1\times S^1\times S^1)$, where $S^1\times S^1\times S^1$ is the boundary of the pair of pants.

This work is part of an ongoing ef\/fort to understand the structure of quantum Riemann surfaces and their noncommutative
dif\/ferential geometry, see~\cite{KL, KL1, KL3, KL4, KM1, KM2,KM3,KM4}.
Our paper has many things in common with the work of Abrahamse~\cite{Abr} and Abrahamse--Douglas~\cite{AD}, who use
dif\/ferent Hilbert spaces.

The paper is organized as follows.
Section~\ref{section2} contains an overview of the zero and one-hole cases, while Section~\ref{section3} has a~detailed
discussion of the quantum pair of pants.

\section{Preliminaries}\label{section2}
In this section we describe in some detail, the zero and one-hole cases.
Most of the material is well-known, however the treatment of the quantum annulus is somewhat new.

\subsection{The quantum disk}

In this subsection we look at the structure of the quantum disk.
We review the tools and the relevant theorems that will be a~motivation for the subsequent discussion of the quantum
pair of pants.

Consider the closed unit disk ${\mathbb D}=\{\zeta\in{\mathbb C}: |\zeta|\leq 1\}$.
We can represent any holomorphic function inside the disk as a~convergent power series
\begin{gather*}
f(\zeta)=\sum\limits_{n=0}^\infty e_n\zeta^n.
\end{gather*}
The Hardy space on the disk is def\/ined as
\begin{gather*}
H^2({\mathbb D})=\left\{f(\zeta)=\sum\limits_{n=0}^\infty e_n\zeta^n: \sum\limits_{n=0}^\infty|e_n|^2
<\infty\right\}.
\end{gather*}
We def\/ine the multiplication operator by the complex coordinate, $z:H^2({\mathbb D}) \to H^2({\mathbb D})$ by the
formula $f(\zeta)\mapsto \zeta f(\zeta)$.
If $E_n=\zeta^n$ is the orthonormal basis on $H^2({\mathbb D})$, then applying~$z$ to the basis elements produces
$zE_n=E_{n+1}$ for all $n \geq 0$, i.e.,~$z$ is the unilateral shift; moreover, we have the following formula for the adjoint operator to~$z$
\begin{gather*}
z^*E_n=
\begin{cases}
E_{n-1} &\text{for} \quad  n\ge 1,
\\
0  &\text{for} \quad n=0.
\end{cases}
\end{gather*}

Now we consider the $C^*$-algebra generated by~$z$.
This well-known algebra is called the Toeplitz algebra, denoted by $\mathcal{T}$, and has also been termed the quantum
(noncommutative) disk.
This is (partially) based on the following standard results collected here with sketches of proofs which serve as
a~guideline for considerations in the next section.

\begin{Theorem}
The norm of~$z$ is~$1$.
The spectrum of~$z$ is all of ${\mathbb D}$, i.e.,~$\sigma(z)={\mathbb D}$.
\end{Theorem}

\begin{proof}
The norm computation is straightforward.
By the norm calculation it then follows that the spectrum is a~closed subset of the unit disk.
To illustrate that any~$\lambda$ in the interior of ${\mathbb D}$ is an eigenvalue of $z^*$, take $f_\lambda(\zeta) =
\sum\limits_{n=0}^\infty \lambda^n \zeta^n$ and so
\begin{gather*}
z^*f_\lambda(\zeta)=\sum\limits_{n=1}^\infty \lambda^n \zeta^{n-1}=\lambda\sum\limits_{n=1}^\infty
\lambda^{n-1}\zeta^{n-1}=\lambda\sum\limits_{n=0}^\infty \lambda^n\zeta^n=\lambda f_\lambda(\zeta).\tag*{\qed}
\end{gather*}
\renewcommand{\qed}{}
\end{proof}

Let $\mathcal{K}$ be the algebra of compact operators in $H^2({\mathbb D})$.
The next observation tells us how the commutator ideal of $\mathcal{T}$, and $\mathcal{K}$ are related.

\begin{Theorem}
The commutator ideal of $\mathcal{T}$ is the ideal of compact operators.
\end{Theorem}

\begin{proof}
Since $\mathcal{T}$ is generated by~$z$ and $z^*$, the commutator ideal of $\mathcal{T}$ is equal to the ideal generated
by the commutator $[z^*,z]$.
Note that $[z^*,z]=P_{E_0}$, the orthogonal projection onto the span of $E_0$.
Since this one-dimensional projection is a~compact operator, it follows that the commutator ideal of $\mathcal{T}$ is
contained in $\mathcal{K}$.
To prove the opposite inclusion we look at the following rank one operators: $E_{ij}(f)=\langle f, E_i\rangle E_j$.
Notice that $E_{ij}=z^jP_{E_0}(z^*)^i$, hence those operators belong to the commutator ideal of $\mathcal{T}$.
But every compact operator is a~norm limit of f\/inite rank operators, which in turn are f\/inite linear combinations of
$E_{ij}$'s.
This verif\/ies that the commutator ideal of $\mathcal{T}$ contains $\mathcal{K}$.
\end{proof}

In order to state the next result, f\/irst we introduce some more notation.
We identify $H^2({\mathbb D})$, the Hardy space on the unit disk, with the subspace of $L^2(S^1)$ spanned~by
$\{e^{inx}\}_{n\ge0}$.
Also given a~continuous function~$f$ on the unit circle, we denote the multiplication operator by~$f$ as~$M_f$.
Let $P:L^2(S^1)\to H^2({\mathbb D})$ be the orthogonal projection onto $\spn\{e^{inx}\}_{n\ge0}$, then def\/ine
the operator $T_f: H^2({\mathbb D}) \to H^2({\mathbb D})$ by~$T_f=PM_f$.
The operator~$T_f$ is known as a~Toeplitz operator.
Since $\|M_f\|=\|f\|_{\infty}$, $\|T_f\|\le \|f\|_{\infty}$ and hence it is bounded.
We have:

\begin{Theorem}
\label{toeplitz_disk}
The quotient $\mathcal{T}/\mathcal{K}$ is isomorphic to $C(S^1)$, the space of continuous functions on the unit circle.
\end{Theorem}

\begin{proof}
The usual proof constructs an isomorphism between the two algebras.
Notice that for a~continuous function~$f$, we have $T_f\in\mathcal{T}$ and since $T_{e^{ix}}$ is the unilateral shift,
$T_{e^{-ix}}=T_{e^{ix}}^*$.
By the Stone--Weierstrass theorem, every continuous function can be approximated by trigonometric polynomials.
Consequently we can def\/ine a~map $\theta: C(S^1) \to \mathcal{T}/\mathcal{K}$ by $\theta: f \mapsto [T_f]$, the class of
operators $T_f$.

Next we show that $T_f$ is compact if and only if $f\equiv 0$.
Suppose $T_f$ is compact.
Then for a~continuous~$f$ with Fourier series $ \sum\limits_{n=-\infty}^\infty e_n e^{inx}$ we have
\begin{gather*}
T_f\big(e^{ikx}\big)=\sum\limits_{n=0}^\infty e_{n-k}e^{inx}.
\end{gather*}
Thus, the matrix coef\/f\/icients $e_n=(E_{i+n},T_f E_i)$ and since $T_f$ is compact, we must have $(E_{i+n}$,
$T_f E_i)\to0$
as $i\to\infty$ for each f\/ixed~$n$.
Therefore, $e_n=0$ for all~$n$ and hence $f\equiv 0$.
This result means that~$\theta$ is injective.

Next we observe that $T_fT_g-T_{fg}$ is a~compact operator for all continuous $f$, $g$.
If $f$, $g$ are trigonometric polynomials then a~direct calculation shows that $T_fT_g-T_{fg}$ is a~f\/inite rank operator.
The general case then follows by appealing to the Stone--Weierstrass theorem.
As a~consequence, the map~$\theta$ above is a~$C^*$-homomorphism.

The range of~$\theta$ is dense since it contains (the classes of) polynomials in~$z$ and $z^*$.
Then by general $C^*$-algebra theory (see~\cite{Conway} for example)~$\theta$ is an isometry hence the range is closed.
This means that $\Ran(\theta)=\mathcal{T}/\mathcal{K}$ and therefore~$\theta$ is a~$*$-isomorphism.
\end{proof}

Note that from the last theorem we get a~short exact sequence
\begin{gather*}
0\rightarrow \mathcal{K}\rightarrow\mathcal{T}\rightarrow C(S^1)\rightarrow 0.
\end{gather*}
We can compare this to the short exact sequence for the classical disk
\begin{gather*}
0\rightarrow C_0({\mathbb D})\rightarrow C({\mathbb D})\rightarrow C(S^1)\rightarrow 0,
\end{gather*}
where $C_0({\mathbb D})$ are the continuous functions on the disk that vanish on the boundary.

\subsection{The quantum annulus}

Let $0<r<1$ and consider the annulus
\begin{gather*}
A_{r}= \left\{\zeta\in{\mathbb C}: r \leq |\zeta| \leq 1\right\}.
\end{gather*}
The classical uniformization theory of Riemann surfaces implies that every open annulus is biholomorphically equivalent
to an annulus of the above form.

We can write any holomorphic function $\varphi(\zeta)$ on the interior of $A_{r}$ as the following convergent version of
Laurent series
\begin{gather*}
\varphi(\zeta)=\sum\limits_{n=0}^\infty e_n\zeta^n+\sum\limits_{n=-\infty}^{-1} f_n\left(\frac{\zeta}{r}\right)^n.
\end{gather*}
We label the basic monomials in the above expansion as
\begin{gather*}
E_n=\zeta^n,
\qquad
F_n =\left(\frac{\zeta}{r}\right)^n,
\end{gather*}
and def\/ine our specially convenient Hilbert space of holomorphic functions on $A_{r}$ to be
\begin{gather*}
H=\left\{\varphi(\zeta)=\sum\limits_{n=0}^\infty e_nE_n+\sum\limits_{n=-\infty}^{-1}f_nF_n:
\|\varphi\|<\infty\right\},
\end{gather*}
where
\begin{gather*}
\|\varphi\|^2=\sum\limits_{n=0}^\infty |e_n|^2+\sum\limits_{n=-\infty}^{-1}|f_n|^2,
\end{gather*}
so that $\{E_n\}$, $\{F_m\}$ form an orthonormal basis.
The operator $z:H \to H$ is def\/ined by the formula $f(\zeta)\mapsto \zeta f(\zeta)$.
With respect to the above basis, the operator~$z$ is a~rather special weighted bilateral shift.
We have
\begin{alignat*}{3}
& zE_n=E_{n+1} \qquad && \text{for} \quad n\ge 0, &
\\
& zF_n=rF_{n+1} \qquad && \text{for} \quad  n\le -2, &
\\
 & zF_{-1}=rE_0&&&
\end{alignat*}
and
\begin{alignat*}{3}
&  z^*E_n=E_{n-1} \qquad && \text{for} \quad n\ge 1,&
\\
& z^*E_0=rF_{-1},&&&
\\
 & z^*F_n=rF_{n-1} \qquad && \text{for} \quad n\le -1.&
\end{alignat*}

In full analogy with the disk case, the operator~$z$ is a~form of a~noncommutative coordinate for what we call quantum
annulus.
First we look at the spectrum of~$z$.

\begin{Theorem}
The norm of~$z$ is~$1$.
The spectrum of~$z$ is all of~$A_{r}$.
\end{Theorem}

\begin{proof}
The formulas above easily imply that $\|z\|\leq 1$, while the action of~$z$ on $E_n$ shows that it is exactly~1.
It is then straightforward to verify that for~$\lambda$ inside $A_{r}$ the following is an eigenvector of $z^*$
corresponding to the eigenvalue~$\lambda$
\begin{gather*}
\phi_\lambda=\sum\limits_{n=0}^\infty \lambda^n E_n+\sum\limits_{n=-\infty}^{-1}\left(\frac{\lambda}{r}\right)^nF_n.
\end{gather*}
Finally, using the techniques described in Lemma~\ref{resolvent_of_z} below, we can prove that the operator $z-\lambda$
is invertible for $|\lambda|<r$.
Put together those statements imply that the spectrum of~$z$ is~$A_{r}$.
\end{proof}

The operators $z^*z$ and $zz^*$ are diagonal.
We have
\begin{alignat*}{3}
 & zz^*E_n=E_{n} \qquad && \text{for} \quad n\ge 1,&
\\
& zz^*F_n=r^2F_{n} \qquad && \text{for}\quad n\le -1,&
\\
& zz^*E_{0}=r^2E_0 &&&
\end{alignat*}
and
\begin{alignat*}{3}
& z^*zE_n=E_{n} \qquad && \text{for} \quad n\ge 1,&
\\
 & z^*zF_n=r^2F_{n} \qquad&& \text{for} \quad n\le -1,&
\\
& z^*zE_0=E_{0}.&&&
\end{alignat*}

Thus the spectrum of those operators is $\sigma(zz^*)=\{1\}\cup\{r^2\}=\sigma(z^*z)$.
Also notice that the spectral projections $P_{z^*z}(1)$ and $P_{z^*z}(r^2)$ of $z^*z$ are orthogonal projections onto
subspaces of~$H$ generated by $E_n$'s and $F_n$'s, respectively.
By the continuous functional calculus applied to $z^*z$, both projections belong to $C^*(z)$, the $C^*$-algebra
generated by~$z$.

\begin{Remark}
The above formulas also imply that the commutator $\frac{z^*z-zz^*}{1-r^2}=\frac{[z^*,z]}{1-r^2}$ is the ortho\-go\-nal
projection onto the one-dimensional subspace spanned by $E_0$,hence a~compact operator.
\end{Remark}

\begin{Theorem}
The commutator ideal of $C^*(z)$ is the ideal of compact operators.
\end{Theorem}
\begin{proof}
By the remark above the commutator ideal of $C^*(z)$ is contained in $\mathcal{K}$.
Similar to the quantum disk case, the opposite inclusion follows from the easily verif\/iable fact that the rank one
operators $f\mapsto\langle f, E_i\rangle E_j$, $f\mapsto\langle f, E_i\rangle F_j$, $f\mapsto\langle f, F_i\rangle E_j$,
$f\mapsto\langle f, F_i\rangle F_j$ are in the commutator ideal of $C^*(z)$.
\end{proof}

\begin{Theorem}
The quotient $C^*(z)/\mathcal{K}$ is isomorphic to $C(S^1)\oplus C(S^1)$, where $C(S^1)$ is the space of continuous
functions on the unit circle.
Thus we have a~short exact sequence
\begin{gather*}
0\rightarrow \mathcal{K}\rightarrow C^*(z)\rightarrow C(S^1)\oplus C(S^1)\rightarrow 0.
\end{gather*}
\end{Theorem}

\begin{proof}
For details we refer to the proof of Theorem~\ref{quotient_algebra} in the next section.
The key step is showing that the inf\/inite-dimensional spectral projections $P_{z^*z}(1)$ and $P_{z^*z}(r^2)$ are in
$C^*(z)$.
They can be used together with Toeplitz operators on subspaces generated by $E_n$'s and $F_n$'s to construct an
isomorphism between $C^*(z)/\mathcal{K}$ and $C(S^1)\oplus C(S^1)$ in a~similar fashion to the Toeplitz algebra case.
\end{proof}

\section{The quantum pair of pants}\label{section3}

Let $0<a<1$, $a+r_2<1$, $r_1+r_2<a$.
We def\/ine the (closed) pair of pants as follows
\begin{gather*}
PP_{(a, r_1,r_2)}= \left\{\zeta\in{\mathbb C}: |\zeta| \leq 1,\, |\zeta|\geq r_1,\, |\zeta-a|\geq r_2\right\}.
\end{gather*}
It is clear that every open disk with two nonintersecting circular holes is biholomorphically equivalent to the interior
of the one of the above pair of pants.
There are some technical advantages to having the holes located as above.
To a~pair of pants we associate a~convenient Hilbert space of holomorphic functions on it and study the operator of
multiplication by~$\zeta$ on that Hilbert space.
This is described more precisely in the following subsection.

\subsection{Definitions}

It follows from~\cite{M} that every holomorphic function on the interior of $PP_{(a, r_1,r_2)}$ can be appro\-xi\-ma\-ted~by
rational functions with the only singularities at the centers of the smaller circles in $PP_{(a, r_1,r_2)}$ or at
inf\/inity.
In fact we can do a~little better.

\begin{Proposition}
Every holomorphic function $\varphi(\zeta)$ on the interior of $PP_{(a, r_1,r_2)}$ can be written as the following
convergent series
\begin{gather*}
\varphi(\zeta)=\sum\limits_{n=0}^\infty e_n\zeta^n+\sum\limits_{n=-\infty}^{-1} f_n\left(\frac{\zeta}{r_1}\right)^n
+ \sum\limits_{n=-\infty}^{-1} g_n\left(\frac{\zeta-a}{r_2}\right)^n.
\end{gather*}
\end{Proposition}

\begin{proof}
In~\cite{Ahlfors}, it was shown that if $\varphi(\zeta)$ is holomorphic on an annulus $\{\zeta\in{\mathbb C}: R_1 <
|\zeta-c|<R_2\}$, then $\varphi(\zeta)=\varphi_1(\zeta)+\varphi_2(\zeta)$ where $\varphi_1(\zeta)$ is holomorphic on
$|\zeta-c|>R_1$ and $\varphi_2(\zeta)$ is holomorphic on $|\zeta-c|<R_2$.
We apply this theorem twice.
Let~$\varphi$ be a~holomorphic function on the open pair of pants.
Consider an annulus $A=\{\zeta\in{\mathbb C}: |\zeta-a|>r_2, |\zeta-c|<r\}$ around~$a$ with outer radius~$r$ and inner
radius~$r_2$ that does not intersect the hole around the origin with radius~$r_1$ and let~$D$ be the disk with
center~$a$ and radius~$r$.
Then $\varphi|_{A}$ is a~holomorphic function and so from~\cite{Ahlfors}, $\varphi|_{A}=\varphi_1+\varphi_2$ with
$\varphi_1$ holomorphic outside the hole centered at~$a$ with radius~$r_2$, and~$\varphi_2$ holomorphic on~$D$.
Consequently~$\varphi_1$ has the following convergent series representation
\begin{gather*}
\varphi_1=\sum\limits_{n=-\infty}^{-1} g_n\left(\frac{\zeta-a}{r_2}\right)^n.
\end{gather*}
Next consider the function $\varphi-\varphi_1$.
This function is holomorphic on $PP_{(a, r_1,r_2)}$ and, because $\varphi-\varphi_1=\varphi_2$ on A, it extends to
a~holomorphic function on~$D$.
This means that $\varphi-\varphi_1$ is holomorphic on the annulus $\{\zeta\in{\mathbb C}: |\zeta|>r_1, |\zeta-c|<1\}$,
and so by using~\cite{Ahlfors} again, we have $\varphi-\varphi_1=\varphi_3+\varphi_4$ with $\varphi_3$ holomorphic
in the unit disk ${\mathbb D}$ and $\varphi_4$ holomorphic on $\{\zeta\in{\mathbb C}: |\zeta|>r_1 \}$.
Thus $\varphi_3$ and $\varphi_4$ have the following convergent series representation
\begin{gather*}
\varphi_3=\sum\limits_{n=0}^\infty e_n \zeta^n
\qquad \text{and} \qquad
\varphi_4=\sum\limits_{n=-\infty}^{-1}f_n\left(\frac{\zeta}{r_1}\right)^n.
\end{gather*}
Combining these three series representations gives the desired result.
\end{proof}

Similar to the annulus case we set
\begin{gather*}
E_n=\zeta^n,
\qquad
F_n =\left(\frac{\zeta}{r_1}\right)^n,
\qquad \text{and} \qquad
G_n=\left(\frac{\zeta-a}{r_2}\right)^n.
\end{gather*}
The Hilbert space~$H$ that we will use is def\/ined as
\begin{gather}
\label{Hilbert_space}
H=\left\{\varphi(\zeta)=\sum\limits_{n=0}^\infty e_nE_n+\sum\limits_{n=-\infty}^{-1}(f_nF_n+g_nG_n):
\|\varphi\|<\infty\right\},
\end{gather}
where
\begin{gather*}
\|\varphi\|^2=\sum\limits_{n=0}^\infty |e_n|^2+\sum\limits_{n=-\infty}^{-1}\left(|f_n|^2+|g_n|^2\right).
\end{gather*}

The advantage of working with the above Hilbert space of holomorphic functions on $PP_{(a, r_1,r_2)}$ is that there is
a~distinguished orthonormal basis in it, namely the basis consisting of $\{E_n\}$, $\{F_m\}$, $\{G_k\}$.

The object of study in this section is the operator $z:H \to H$ given by $z\varphi(\zeta)=M_\zeta \varphi(\zeta) =
\zeta \varphi(\zeta)$, i.e., the multiplication operator by~$\zeta$.
Straightforward calculations yields the following formulas.

\begin{Lemma}
The operators~$z$ and $z^*$ act on the basis elements in the following way
\begin{alignat*}{3}
& zE_n=E_{n+1} \qquad && \text{for} \quad n\ge 0,&
\\
& zF_n=r_1F_{n+1} \qquad && \text{for} \quad n\le -2, &
\\
& zF_{-1}=r_1E_0,&&&
\\
& zG_n=r_2G_{n+1}+aG_n \qquad && \text{for} \quad  n\le -2, &
\\
& zG_{-1}=r_2E_0+aG_{-1}&&&
\end{alignat*}
and
\begin{alignat*}{3}
& z^*E_n=E_{n-1} \qquad && \text{for} \quad  n\ge 1,&
\\
& z^*E_0=r_1F_{-1}+r_2G_{-1}, &&&
\\
& z^*F_n=r_1F_{n-1} \qquad && \text{for} \quad n\le -1,&
\\
& z^*G_n=r_2G_{n-1}+aG_n \qquad && \text{for} \quad n\le -1. &
\end{alignat*}
\end{Lemma}

\begin{Lemma}
\label{z_z_star_coefficients}
The operators~$z$ and $z^*$ shift the coefficients of $\varphi(\zeta)$ in the series decomposition defined in
equation~\eqref{Hilbert_space} in the following way
\begin{gather*}
z\varphi=\sum\limits_{n=0}^\infty\tilde{e}_nE_n+\sum\limits_{n=-\infty}^{-1}\!\! \big(\tilde{f}_nF_n+\tilde{g}_nG_n\big)
\qquad
\text{and}
\qquad
z^*\varphi=\sum\limits_{n=0}^\infty e'_nE_n+\sum\limits_{n=-\infty}^{-1}\!\! (f'_nF_n+g'_nG_n),
\end{gather*}
where
\begin{alignat*}{3}
& \tilde{e}_n=e_{n-1} \qquad && \text{for} \quad n\ge 1, &
\\
& \tilde{e}_0=r_1f_{-1}+r_2g_{-1}, &&&
\\
& \tilde{f}_n=r_1f_{n-1} \qquad && \text{for} \quad n\le -1,&
\\
& \tilde{g}_n=r_2g_{n-1}+ag_n \qquad && \text{for} \quad n\le -1 &
\end{alignat*}
and
\begin{alignat*}{3}
& e'_n=e_{n+1} \qquad && \text{for} \quad n\ge 0, &
\\
& f'_n=r_1f_{n+1} \qquad && \text{for} \quad n\le -2, &
\\
& f'_{-1}=r_1e_0, &&&
\\
& g'_n=r_2g_{n+1}+ag_n \qquad && \text{for} \quad n\le -2,&
\\
& g'_{-1}=r_2e_0+ag_{-1}. &&&
\end{alignat*}
\end{Lemma}
We can now def\/ine the quantum pair of pants.

\begin{Definition}
The quantum pair of pants, denoted $QPP_{(a, r_1,r_2)}$, is def\/ined to be the $C^*$-algebra generated by the
operator~$z$,
i.e., $QPP_{(a, r_1,r_2)}=C^*(z)$.
\end{Definition}

\subsection[The spectrum of~$z$]{The spectrum of~$\boldsymbol{z}$}

In this subsection we study the spectrum of~$z$, starting with a~calculation of the norm of~$z$.

\begin{Proposition}
\label{norm_of_z}
With the above notation, we have: $\|z\|=1$.
\end{Proposition}

\begin{proof}
Using the series representation of $\varphi(\zeta)$ in formula \eqref{Hilbert_space} above, and the coef\/f\/icients of
Lemma~\ref{z_z_star_coefficients} we compute $\|z\varphi\|^2$:
\begin{gather*}
\|z\varphi\|^2=\sum\limits_{n=1}^\infty |e_{n-1}|^2+|r_1f_{-1}+r_2g_{-1}|^2+r_1^2\sum\limits_{n=-\infty}^{-1}
|f_{n-1}|^2+\sum\limits_{n=-\infty}^{-1}|r_2g_{n-1}+ag_n|^2.
\end{gather*}
Using the triangle inequality and the fact that $a> r_1$ we obtain
\begin{gather*}
\|z\varphi\|^2 \le \sum\limits_{n=0}^\infty |e_n|^2+\big(a|f_{-1}|+r_2|g_{-1}|\big)^2+r_1^2
\sum\limits_{n=-\infty}^{-1}|f_{n-1}|^2+\sum\limits_{n=-\infty}^{-1}\big(r_2|g_{n-1}|+a|g_n|\big)^2.
\end{gather*}
Notice that by denoting $g_0:= f_{-1}$ we can write
\begin{gather*}
\big(a|f_{-1}|+r_2|g_{-1}|\big)^2+\sum\limits_{n=-\infty}^{-1}\big(r_2|g_{n-1}|+a|g_n|\big)^2=
\sum\limits_{n=-\infty}^{0}\big(r_2|g_{n-1}|+a|g_n|\big)^2
\\
\qquad{}
=r_2^2 \sum\limits_{n=-\infty}^0 |g_{n-1}|^2+ a^2\sum\limits_{n=-\infty}^0 |g_n|^2+2ar_2 \sum\limits_{n=-\infty}^0
|g_{n-1}||g_n|.
\end{gather*}
The Cauchy--Schwartz inequality implies
\begin{gather*}
\sum\limits_{n=-\infty}^{0}\big(r_2|g_{n-1}|+a|g_n|\big)^2   \leq r_2^2 \sum\limits_{n=-\infty}^0 |g_{n-1}|^2+
a^2\sum\limits_{n=-\infty}^0 |g_n|^2
\\
\qquad
\phantom{\leq}  {}
+2ar_2 \left(\sum\limits_{n=-\infty}^0 |g_{n-1}|^2\right)^{1/2}\left(\sum\limits_{n=-\infty}^0 |g_n|^2\right)^{1/2}
\\
\qquad
\leq \big(r_2^2+a^2+2ar_2\big) \sum\limits_{n=-\infty}^0 |g_n|^2 =(r_2+a)^2 \left(|f_{-1}|^2+\sum\limits_{n=-\infty}^{-1} |g_n|^2\right).
\end{gather*}
Using the fact that $r_1, r_2+a <1$ in the above computations we see that
\begin{gather*}
\|z\varphi\|^2  \le \sum\limits_{n=0}^\infty |e_n|^2+r_1^2\sum\limits_{n=\infty}^{-1}|f_{n-1}|^2+(r_2+a)^2
|f_{-1}|^2+(r_2+a)^2\sum\limits_{n=-\infty}^{-1}|g_n|^2
\\
\phantom{\|z\varphi\|^2}{}
 \leq \sum\limits_{n=0}^\infty |e_n|^2+\sum\limits_{n=\infty}^{-1}|f_n|^2+\sum\limits_{n=-\infty}^{-1}|g_n|^2 =
\|\varphi\|^2,
\end{gather*}
showing that $\|z\|\le 1$.
On the other hand, $\|zE_1\| =\| E_2\|=\|E_1\|$.
Thus $\|z\|=1$.
\end{proof}

Next we compute the spectrum of~$z$.
In estimating the norms of resolvents of~$z$ we use the following well known result.

\begin{Lemma}[Schur--Young inequality]
\label{schuryounginq}
Let $T: L^2(Y) \longrightarrow L^2(X)$ be an integral operator
\begin{gather*}
Tf(x)=\int K(x,y)f(y)dy.
\end{gather*}
Then one has
\begin{gather*}
\|T\|^2 \le \left(\sup_{x\in X} \int_Y |K(x,y)|dy\right)\left(\sup_{y\in Y} \int_X |K(x,y)|dx\right).
\end{gather*}
\end{Lemma}

The details of the lemma and its proof can be found in~\cite{HS}.

\begin{Lemma}
\label{resolvent_of_z}
The operator $z-\lambda$ has a~bounded inverse for $|\lambda|< r_1$, $|\lambda-a|< r_2$, and $|\lambda|>1$.
\end{Lemma}

\begin{proof}
Let
\begin{gather*}
\varphi(\zeta)=\sum\limits_{n=0}^\infty e_n E_n+\sum\limits_{n=-\infty}^{-1}(f_nF_n+g_nG_n)
\end{gather*}
and
\begin{gather*}
\tilde{\varphi}(\zeta)=\sum\limits_{n=0}^\infty \tilde{e}_n E_n+\sum\limits_{n=-\infty}^{-1}\big(\tilde{f}_nF_n+\tilde{g}_nG_n\big).
\end{gather*}
Consider the equation $(z-\lambda)\varphi(\zeta)=\tilde{\varphi}(\zeta)$.
Using the above decompositions and Lemma~\ref{z_z_star_coefficients} we obtain the following system of equations
\begin{alignat}{3}
&  r_1f_{-1}+r_2g_{-1}-\lambda e_0=\tilde{e}_0,\qquad &&&\nonumber
\\
& e_{n-1}-\lambda e_n=\tilde{e}_n \quad &&\text{for} \quad n\ge1, &\nonumber
\\
& r_1f_{n-1}-\lambda f_n=\tilde{f}_n \quad &&\text{for} \quad n\le-1,&\nonumber
\\
 & r_2g_{n-1}+ag_n-\lambda g_n=\tilde{g}_n \quad &&\text{for} \quad  n\le-1. & \label{resolvent_system}
\end{alignat}

By Proposition~\ref{norm_of_z}, $\|z\|=1$ and if $|\lambda|>1=\|z\|$ then by general functional analysis we know that
$(z-\lambda)^{-1}$ is a~bounded, invertible operator.

Next we consider three cases: the f\/irst case is for $0<|\lambda|<r_1$, the second case is for $|\lambda-a|<r_2$, and the
last case is for $\lambda=0$.

If $0<|\lambda|<r_1<1$, then $|\lambda-a|>r_2$.
We can solve the system of equations~\eqref{resolvent_system} recursively.
Rewriting the last equation and multiplying by $((\lambda-a)/r_2)^{n-1}$ yields
\begin{gather*}
\left(\frac{\lambda-a}{r_2}\right)^{n-1}g_{n-1}-\left(\frac{\lambda-a}{r_2}\right)^ng_n
=\left(\frac{\lambda-a}{r_2}\right)^{n-1}\frac{1}{r_2}\tilde{g}_n.
\end{gather*}
Letting $h_n=((\lambda-a)/r_2)^ng_n$, we get
\begin{gather*}
h_{n-1}-h_n=\left(\frac{\lambda-a}{r_2}\right)^{n-1}\frac{1}{r_2}\tilde{g}_n.
\end{gather*}
The requirement for a~square summable solution forces $h_n =
-\sum\limits_{j=-\infty}^n((\lambda-a)/r_2)^{j-1}\tilde{g}_j/r_2$ and hence for $n\le-1$ we obtain
\begin{gather*}
g_n=-\frac{1}{r_2}\sum\limits_{j=-\infty}^n\left(\frac{\lambda-a}{r_2}\right)^{j-n-1}\tilde{g}_j.
\end{gather*}
Similar calculations show that
\begin{gather*}
e_n=\sum\limits_{j=n+1}^\infty \lambda^{j-n-1}\tilde{e}_j
\qquad
\text{and}
\qquad
f_n=\left(\frac{\lambda}{r_1}\right)^{-n-1}f_{-1}+\frac{1}{r_1}\sum\limits_{j=n+1}^{-1}\left(\frac{\lambda}{r_1}\right)^{j-n-1}\tilde{f}_j
\end{gather*}
for $n\ge0$ and $n\le-2$ respectively.
These formulas along with the f\/irst equation in system~\eqref{resolvent_system} give
\begin{gather*}
f_{-1}=\frac{1}{r_1}\left(\sum\limits_{j=0}^\infty\lambda^j\tilde{e}_j+\sum\limits_{j=-\infty}^{-1}\left(\frac{\lambda
- a}{r_2}\right)^j\tilde{g}_j\right).
\end{gather*}

We introduce some notation; f\/irst notice that we have a~natural decomposition, $H \!\cong\! \ell^2({\mathbb Z}_{\ge0})
\oplus \ell^2({\mathbb Z}_{<0}) \oplus \ell^2({\mathbb Z}_{<0})$ given in the following way: for $\varphi\in H$ write
$\varphi=e+f+g$ where $e=\sum\limits_{n\ge0}e_nE_n$, $f= \sum\limits_{n\le-1}f_nF_n$, and $g=\sum\limits_{n\le-1}g_nG_n$.
Using this notation we see that $\|\varphi\|^2=\|e\|^2+\|f\|^2+\|g\|^2$.
Def\/ine the characteristic $\chi(t)=1$ for $0\le t\le1$ and zero otherwise, then we can def\/ine seven dif\/ferent integral operators
\begin{gather}
T_1e=\sum\limits_{n=0}^\infty\sum\limits_{j=0}^\infty \lambda^{j-n-1}\chi\left(\frac{n+1}{j}\right)e_jE_n: \
\ell^2({\mathbb Z}_{\ge0})\to\ell^2({\mathbb Z}_{\ge0}),
\nonumber
\\
T_2(e,g)=\sum\limits_{n=-\infty}^{-1}\frac{1}{r_1}\left(\frac{\lambda}{r_1}\right)^{-n-1}\left(\sum\limits_{j=0}^\infty\lambda^je_j
+ \sum\limits_{j=-\infty}^{-1}\left(\frac{\lambda-a}{r_2}\right)^jg_j\right)F_n :
\nonumber
\\
\hphantom{T_2(e,g)=}{} \ell^2({\mathbb Z}_{\ge0})\oplus\ell^2({\mathbb Z}_{<0})\to\ell^2({\mathbb Z}_{<0})
\nonumber
\\
T_3f =
\sum\limits_{n=-\infty}^{-1}\frac{1}{r_1}\sum\limits_{j=-\infty}^{-1}\left(\frac{\lambda}{r_1}\right)^{j-n-1}\chi\left(\frac{n+1}{j}\right)f_jF_n: \
\ell^2({\mathbb Z}_{<0})\to\ell^2({\mathbb Z}_{<0}),
\nonumber
\\
\label{int_ops_resolvent_z}
T_4f =
\sum\limits_{n=-\infty}^{-1}\frac{1}{r_1}\sum\limits_{j=-\infty}^{-1}\left(\frac{\lambda}{r_1}\right)^{j-n-1}\chi\left(\frac{j}{n}\right)f_jF_n: \
\ell^2({\mathbb Z}_{<0})\to\ell^2({\mathbb Z}_{<0}),
\\
T_5(e,f) =\sum\limits_{n=-\infty}^{-1}\frac{1}{r_2}\left(\frac{\lambda-a}{r_2}\right)^{-n-1}
\left(\sum\limits_{j=0}^\infty\lambda^je_j+\sum\limits_{j=-\infty}^{-1}\left(\frac{\lambda}{r_1}\right)^jf_j\right)G_n:
\nonumber
\\
\hphantom{T_5(e,f) =}{} \ell^2({\mathbb Z}_{\ge0})\oplus\ell^2({\mathbb Z}_{<0})\to\ell^2({\mathbb Z}_{<0}),
\nonumber
\\
T_6g=\sum\limits_{n=-\infty}^{-1}\frac{1}{r_2}\sum\limits_{j=-\infty}^{-1}\left(\frac{\lambda-a}{r_2}\right)^{j-n-1}
\chi\left(\frac{n+1}{j}\right)g_j G_n: \ \ell^2({\mathbb Z}_{<0})\to\ell^2({\mathbb Z}_{<0}),
\nonumber
\\
T_7g=\sum\limits_{n=-\infty}^{-1}\frac{1}{r_2}\sum\limits_{j=-\infty}^{-1}\left(\frac{\lambda-a}{r_2}\right)^{j-n-1}
\chi\left(\frac{j}{n}\right)g_j G_n: \  \ell^2({\mathbb Z}_{<0})\to\ell^2({\mathbb Z}_{<0}).
\nonumber
\end{gather}

The operators from formula~\eqref{int_ops_resolvent_z} can be used to represent $(z-\lambda)^{-1}\tilde{\varphi}$, for
$\tilde{\varphi}=\tilde{e}+\tilde{f}+\tilde{g}$ where $\tilde{e}=\sum\limits_{n\ge0}\tilde{e}_nE_n$,
$\tilde{f}=\sum\limits_{n\le -1}\tilde{f}_nF_n$ and $\tilde{g}=\sum\limits_{n\le -1}\tilde{g}_nG_n$, in the following way
\begin{gather*}
(z-\lambda)^{-1}\tilde{\varphi}=T_1\tilde{e}+T_2(\tilde{e},\tilde{g})+T_3\tilde{f}-T_7\tilde{g}.
\end{gather*}
Next we estimate the norm of $(z-\lambda)^{-1}$.
We use Lemma~\ref{schuryounginq} to estimate the norms of the operators $T_1$, $T_3$ and $T_7$ and we directly estimate
$\|T_2\tilde{f}\|$.
The f\/irst estimate is
\begin{gather*}
\|T_1\|^2
\leq \left(\sup_{n \geq 0} |\lambda|^{-n-1} \sum\limits_{j=n+1}^{\infty}|\lambda|^j \right)\left(\sup_{j \geq 1} |\lambda|^{j-1}
\sum\limits_{n=0}^{j-1}|\lambda|^{-n} \right)
\\
\phantom{\|T_1\|^2}{}
=\frac 1{\left(1-|\lambda|\right)}\left(\sup_{j \geq 1} \frac{1-|\lambda|^j}{1-|\lambda|} \right)=\frac1{\left(1-|\lambda|\right)^2},
\end{gather*}
where we have used the fact that $|\lambda|<1$.
Similarly, we have
\begin{gather*}
\|T_3\|^2
\leq \frac 1{r_1^2}\left(\sup_{n \leq -2} \left(\frac{|\lambda|}{r_1}\right)^{-n-1}
\sum\limits_{j=n+1}^{-1}\left(\frac{|\lambda|}{r_1}\right)^j \right) \left(\sup_{j \leq -1}
\left(\frac{|\lambda|}{r_1}\right)^{j-1} \sum\limits_{n=-\infty}^{j-1}\left(\frac{|\lambda|}{r_1}\right)^{-n} \right)
\\
\phantom{\|T_3\|^2}{}
  \leq \frac 1{r_1^2 \big(1- \frac{|\lambda|}{r_1}\big)^2} \left(\sup_{n\leq -2}
1-\left(\frac{|\lambda|}{r_1}\right)^{-n-1}\right).
\end{gather*}
Since $\frac{|\lambda|}{r_1} <1$, it follows that $\|T_3\|^2 \leq \frac{1}{r_1^2\big(1- \frac{|\lambda|}{r_1}\big)^2}$.

Next,
\begin{gather*}
\|T_7\|^2  \leq \frac 1{r_2^2} \left(\sup_{n \leq -1} \left|\frac{\lambda -a}{r_2} \right|^{-n-1}
\sum\limits_{j=-\infty}^n \left|\frac{\lambda -a}{r_2} \right|^j \right) \left(\sup_{j \leq -1} \left|\frac{\lambda
-a}{r_2} \right|^{j-1} \sum\limits_{n=j}^{-1} \left|\frac{\lambda -a}{r_2} \right|^{-n} \right)
\\
\phantom{\|T_7\|^2}{}
  \leq \frac 1{r_2^2} \frac{(|\lambda -a|/ r_2)^{-2}}{\big(1- (|\lambda -a|/ r_2)^{-1}\big)^2} \left(\sup_{j \leq -1} 1-
\left(\frac{|\lambda -a|}{r_2}\right)^j\right).
\end{gather*}
Because $\frac{|\lambda -a|}{r_2}>1$, we have
\begin{gather*}
\|T_7\|^2 \leq \frac 1{r_2^2} \frac{(|\lambda -a|/ r_2)^{-2}}{\big(1- (|\lambda -a|/ r_2)^{-1}\big)^2}=\frac {1}{r_1^2
\big(\frac{|\lambda -a|}{r_2} -1\big)^2}.
\end{gather*}

The operator $T_2$ is a~rank one operator and the norm $T_2\tilde{f}$ can be estimated directly, using the
Cauchy--Schwartz inequality
\begin{gather*}
\|T_2(\tilde{e},\tilde{g})\|^2=
\frac{1}{r_1^2}\sum\limits_{n=-\infty}^{-1}\left(\frac{|\lambda|}{r_1}\right)^{-2n-2}\left|\sum\limits_{j=0}^\infty\lambda^j\tilde{e}_j
+ \sum\limits_{j=-\infty}^{-1}\left(\frac{\lambda-a}{r_2}\right)^j\tilde{g}_j\right|^2
\\
\phantom{\|T_2(\tilde{e},\tilde{g})\|^2}{}
 \le\frac{1}{r_1^2(1-(|\lambda|/r_1)^2)}\left(\frac{1}{1-|\lambda|^2}+ \frac{1}{(|\lambda-a|/r_2)^2-1}\right)
\left(\|\tilde{e}\|^2+\|\tilde{g}\|^2\right).
\end{gather*}
This shows that $(z-\lambda)^{-1}$ is bounded for $0<|\lambda|<r_1$.

The second case is $|\lambda-a|<r_2$.
This implies that $r_1<|\lambda|<1$.
Under these constraints we solve system~\eqref{resolvent_system} using the same methods as those for the f\/irst case to
obtain
\begin{gather*}
e_n =\sum\limits_{j=n+1}^\infty \lambda^{j-n-1}\tilde{e}_j
\qquad\text{for}\quad n\ge0,
\\
f_n =-\frac{1}{r_1}\sum\limits_{j=-\infty}^n \left(\frac{\lambda}{r_1}\right)^{j-n-1}\tilde{f}_j
\qquad\text{for}\quad n\le-1,
\\
g_n =-\left(\frac{\lambda-a}{r_2}\right)^{-n-1}g_{-1}+\frac{1}{r_2}\sum\limits_{j=n+1}^{-1}\left(\frac{\lambda-a}{r_2}\right)^{j-n-1}\tilde{g}_j
\qquad
\text{for}
\quad
n\le-2.
\end{gather*}
Then the f\/irst equation of system~\eqref{resolvent_system} gives
\begin{gather*}
g_{-1}=\frac{1}{r_2}\left(\sum\limits_{j=0}^\infty\lambda^j\tilde{e}_j+\sum\limits_{j=-\infty}^{-1}\left(\frac{\lambda}{r_1}\right)^j\tilde{f}_j\right).
\end{gather*}

Similar to the f\/irst case we can express $(z-\lambda)^{-1}$ using the operators def\/ined in
formula~\eqref{int_ops_resolvent_z} to get
\begin{gather*}
(z-\lambda)^{-1}\tilde{\varphi}=T_1\tilde{e}-T_4\tilde{f}+T_5(\tilde{e},\tilde{f})+T_6\tilde{g}.
\end{gather*}

We omit the repetitive details of estimates of $T_4$, $T_5$, and $T_6$ norms.
They imply that $(z-\lambda)^{-1}$ is bounded for $|\lambda -a|<r_2$.

The last case is when $\lambda=0$.
Solving system~\eqref{resolvent_system} we obtain
\begin{gather*}
e_n =\tilde{e}_{n+1} \qquad\text{for}\quad n\ge0,
\\
f_n =\frac{1}{r_1}\tilde{f}_{n+1} \qquad\text{for}\quad n\le-2,
\\
g_n =-\frac{1}{r_2}\sum\limits_{j=-\infty}^n\left(-\frac{a}{r_2}\right)^{j-n-1}\tilde{g}_j \qquad\text{for}\quad n\le-1.
\end{gather*}
Using the f\/irst equation of system~\eqref{resolvent_system} we compute $f_{-1}$
\begin{gather*}
f_{-1}=\frac{1}{r_1}\left(\tilde{e}_0+\sum\limits_{j=-\infty}^{-1}\left(-\frac{a}{r_2}\right)^j\tilde{g}_j\right).
\end{gather*}
As before the norm estimates hinge on convergent geometric series.
This completes the proof.
\end{proof}

\begin{Theorem}
The spectrum of~$z$ is the regular pair of pants, i.e., $\sigma(z)=PP_{(a, r_1,r_2)}$.
\end{Theorem}

\begin{proof}
By Proposition~\ref{norm_of_z}, $\sigma(z)\subset {\mathbb D}$.
Let
\begin{gather*}
\varphi_\lambda(\zeta)=\sum\limits_{n=0}^\infty \lambda^nE_n+\sum\limits_{n=-\infty}^{-1}
\left(\frac{\lambda}{r_1}\right)^nF_n+\sum\limits_{n=-\infty}^{-1}\left(\frac{\lambda-a}{r_2}\right)^nG_n.
\end{gather*}
It is easy to see that for any~$\lambda$ in the interior of $PP_{(a, r_1,r_2)}$, $\varphi_\lambda\in H$ and
that~$\lambda$ is an eigenvalue with $\varphi_\lambda$ as the associated eigenfunction for $z^*$.
Therefore, $PP_{(a, r_1,r_2)}\subset\sigma(z)$.
From Lemma~\ref{resolvent_of_z} the operator $z-\lambda$ has a~bounded inverse whenever $|\lambda|< r_1$ or $|\lambda-a|<r_2$.
Hence the resolvent set is contained in the holes within the unit disk or outside the unit disk, and so
$\sigma(z)\subset PP_{(a, r_1,r_2)}$.
\end{proof}

In view of the above theorem we can think of the operator~$z$ as a~form of a~noncommutative complex coordinate for what
we call quantum pair of pants.

\subsection[Structure of $C^*(z)$]{Structure of $\boldsymbol{C^*(z)}$}

Next we study the commutator ideal of $C^*(z)$.
A~straightforward computation gives the following formulas.

\begin{Lemma}
\label{commutator_calculation}
The commutator of $z^*$ and~$z$ act on the basis elements in the following way: $[z^*,z]E_n=0$ for $n\ge 1$,
$[z^*,z]F_n=0$ for $n\le -2$, $[z^*,z]G_n =0$ for $n\le -2$.
Moreover on the initial elements we get:
\begin{gather*}
 [z^*,z]E_0=\big(1-r_1^2-r_2^2\big)E_0-ar_2G_{-1},
\\
 [z^*,z]F_{-1} =r_1r_2G_{-1},
\\
 [z^*,z]G_{-1} =r_1r_2F_{-1}-ar_2E_0.
\end{gather*}
\end{Lemma}

Let $\mathcal{I}$ be the ideal generated by $[z^*,z]$.
It is easy to see that $\mathcal{I}$ is in fact the commutator ideal of $C^*(z)$ because that algebra is singly generated.

\begin{Theorem}
The commutator ideal $\mathcal{I}$ of $C^*(z)$ is the $C^*$-algebra $\mathcal{K}$ of compact operators in~$H$.
\end{Theorem}

\begin{proof}
From Lemma~\ref{commutator_calculation} it is clear that the commutator $[z^*, z]$ is f\/inite-rank and hence compact.
Thus, $\mathcal I \subset \mathcal K$.
On the other hand, to show that $\mathcal K \subset \mathcal I$ we will use the following step by step method building
up to the conclusion that a~large collection of rank one operators belong to $\mathcal I$ and that the compact operators
are exactly the norm limit of those.

{\bf Step 1.} First we show that $P=$ orthogonal projection onto $\operatorname{span}\{E_0, F_{-1}, G_{-1}\}$ belongs to the commutator
ideal.
Notice that the (self-adjoint) operator $[z^*{,} z]$ acting on $\operatorname{span}\{E_0, F_{{-}1}, G_{{-}1}\}$ has the following matrix
representation in the basis $\{E_0, F_{-1}, G_{-1}\}$:
\begin{gather*}
A= \left(
\begin{matrix}
1-r_1^2-r_2^2 & 0 & -ar_2
\\
0 & 0 & r_1r_2
\\
-ar_2 & r_1r_2 & 0
\end{matrix}
\right).
\end{gather*}

This matrix has rank equal to $3$ and the following characteristic polynomial
\begin{gather*}
p_A(\lambda)=\lambda^3-\big(1-r_1^2-r_2^2\big)\lambda^2-\big(a^2r_2^2+r_1^2r_2^2\big)\lambda+\big(1-r_1^2-r_2^2\big)r_1^2r_2^2.
\end{gather*}

Since $0 < r_1$, $r_2 < 1$ it is clear that zero is not an eigenvalue of~$A$.
If $\lambda_i$, $i=1,2,3$ are the roots of $p_A(\lambda)$ then by functional calculus there exists a~continuous function
$f: \mathbb R \to \mathbb R$ such that $f(0)=0$, $f(\lambda_i)=1$ so that $f([z^*,z])=P$.
Consequently $P \in \mathcal I \subset C^*(z)$.

{\bf Step 2.} The next step is showing that $P_{E_1}=$ orthogonal projection onto span of $\{E_1\}$ belongs to $\mathcal I$.
We f\/irst observe that the operator $zPz^*$ acts on the basis elements in the following way
\begin{gather*}
zpz^* B=
\begin{cases}
\big(r_1^2+r_2^2\big)E_0+ar_2G_{-1} & \text{if}\quad  B=E_0,
\\
E_1 & \text{if}\quad  B=E_1,
\\
ar_2E_0+a^2 G_{-1} & \text{if}\quad  B=G_{-1},
\\
0 & \text{otherwise}.
\end{cases}
\end{gather*}
Thus, the operator $zPz^*$ on $\operatorname{span}\{E_0, G_{-1}\}$ is self-adjoint and has the following matrix representation in the
basis $\{E_0, G_{-1}\}$:
\begin{gather*}
C=\left(
\begin{matrix}
r_1^2+r_2^2 & ar_2
\\
ar_2 & a^2
\end{matrix}
\right).
\end{gather*}
The characteristic polynomial for~$C$ is $p_C(\lambda)=\lambda^2-(r_1^2+r_2^2+a^2)\lambda+a^2r_1^2$.
First we need to show that $\lambda=0,1$ are not roots of $p_C(\lambda)$.
Clearly $p_C(0) \neq 0$.
Suppose $p_C(1)=0$.
Then solving the equation for $r_2$ we obtain, $r_2^2=(1-a^2)(1-r_1^2)$.
Since $r_2 < 1-a$ and $r_2 < 1-r_1$ we see that
\begin{gather*}
\big(1-a^2\big)\big(1-r_1^2\big) < (1-a)(1-r_1),
\end{gather*}
implying $(1+a)(1+r_1)<1$ which is clearly a~contradiction since $a, r_1 >0$.
Thus, $p_C(1)\neq 0$.

Now we look at the discriminant~$\Delta$ of $p_C(\lambda)$:
\begin{gather*}
\Delta=\big(r_1^2+r_2^2+a^2\big)^2-4a^2r_1^2.
\end{gather*}

If $\Delta=0$ this would imply that $a=r_1$ and $r_2=0$, which is a~contradiction.
Hence~$C$ has two distinct eigenvalues $\lambda_1$ and $\lambda_2$.
Thus, once again by functional calculus there exists a~continuous real valued function~$f$ such that $f(0)=
f(\lambda_1)=f(\lambda_2)=0$ and $f(1)=1$.
Consequently, applying~$f$ to $zPz^*$ we get that, $f(zPz^*)= p_{E_1} \in \mathcal I$.

{\bf Step 3.} By similar functional calculus argument as above we also see that $P_{E_0, G_{-1}}$, the orthogonal
projection onto $\operatorname{span}\{E_0, G_{-1}\}$, belongs to $\mathcal I$.
Consequently, if $P_{F_{-1}}=$ orthogonal projection onto span of $\{F_{-1}\}$ then clearly, $P_{F_{-1}}= P- P_{E_0,
G_{-1}} \in \mathcal I$.

{\bf Step 4.} We will show that $P_{E_n}= $ orthogonal projection onto $\operatorname{span}\{E_n\}$ belongs to $\mathcal I$ for $n=1,2,
\dots$.
To this end, we compute the action of $z^{n-1} P_{E_1} (z^*)^{n-1}$ on the basis elements
\begin{gather*}
z^{n-1} P_{E_1} (z^*)^{n-1} B=
\begin{cases}
E_n & \text{if}\quad  B=E_n,
\\
0 & \text{otherwise.}
\end{cases}
\end{gather*}
Therefore, $z^{n-1} P_{E_1} (z^*)^{n-1}= P_{E_n} \in \mathcal I$ for $n \geq 1$.

{\bf Step 5.} Now consider the action of $(z^*)^n P_{F_{-1}} z^n$, $n \geq 1$ on basis elements
\begin{gather*}
(z^*)^n P_{F_{-1}} z^n B=
\begin{cases}
r_1^{2n}F_{-n-1} & \text{if}\quad  B=F_{-n-1},
\\
0 & \text{otherwise.}
\end{cases}
\end{gather*}
Thus, $r_1^{2n}(z^*)^{-n} P_{F_{-1}} z^{-n}= P_{F_n}$, the projection onto $F_n$, for $n \leq -1$ and hence the
projec\-tions~$P_{F_n}$ belong to $\mathcal I$.

{\bf Step 6.} Since $z^* P_{E_1}z E_0=E_0$ and zero elsewhere, it is clear that $z^* P_{E_1}z =P_{E_0} \in \mathcal I$.
Hence $P_{G_{-1}}= P- P_{E_0}- P_{F_{-1}}$ also belongs to~$\mathcal I$.

{\bf Step 7.} It remains to show that for $n \leq -2$ the orthogonal projection $P_{G_n}$ onto $G_n$ belongs to~$\mathcal I$.
We consider the action of $z^* P_{G_{-1}}z$ on basis elements
\begin{gather*}
z^* P_{G_{-1}}z B=
\begin{cases}
a^2 G_{-1}+ar_2G_{-2} & \text{if}\quad  B=G_{-1},
\\
ar_2 G_{-1}+r_2^2G_{-2} & \text{if}\quad  B= G_{-2},
\\
0 & \text{otherwise}.
\end{cases}
\end{gather*}
Thus it has the following matrix representation relative to the basis $\{G_{-1}, G_{-2}\}$:
\begin{gather*}
D=\left(
\begin{matrix}
a^2 & ar_2
\\
ar_2 & r_2^2
\end{matrix}
\right).
\end{gather*}
The matrix~$D$ has eigenvalues $\lambda_1=0$ and $\lambda_2= a^2+r_2^2$ with $v= aG_{-1}+ r_2 G_{-2}$ being the
eigenvector corresponding to $\lambda_2$.
Thus $P_{v}$, the one-dimensional orthogonal projection onto~$v/\|v\|$, belongs to~$\mathcal I$.
Since $P_{G_{-1}}$ and $P_{v}$ are not mutually orthogonal, simple matrix algebra shows that the set $\{I, P_{G_{-1}},
P_{v}, P_{G_{-1}}P_{v}\}$, where~$I$ is the $2 \times 2$ identity matrix, generates the set of all $2 \times 2 $
matrices.
Consequently, $P_{G_{-2}}$ can be written as a~linear combination of these four projections, making it clear that
$P_{G_{-2}} \in \mathcal I$.

Finally we use induction on~$n$ and follow a~similar argument as above to show that $P_{G_n} \in \mathcal I$ for
$n=-2,-3, -4, \dots$.

{\bf Step 8.} Next we proceed to show that the one-dimensional operators $P_{B_i, B_j}(x)= \langle B_i, x \rangle B_j$,
where $B_i$, $B_j$ are basis elements, i.e., elements of the set $\{E_n, F_k, G_k: n \geq 0, k \leq -1\}$, also belong to
$\mathcal I$.

Since $z^m P_{E_n}E_n= E_{n+m}$ we see that $P_{E_n, E_{n+m}}=z^m P_{E_n}$ for every $n,m \geq 0$.
Similarly we observe that $P_{E_n, E_{n-m}}=(z^*)^mP_{E_n}$ for $m \leq n$.
Together, this proves that all operators $P_{E_n, E_k} \in \mathcal I$ for any $n,k \geq 0$.

Next, we observe that
\begin{gather*}
z^m P_{F_n}=
\begin{cases}
r_1^m F_{n+m} & \text{if}\quad  n+m <0,
\\
r_1^{-n}E_{n+m} & \text{if}\quad  n+m \geq 0.
\end{cases}
\end{gather*}
Moreover, $(z^*)^m P_{F_n}=r_1^m F_{n-m}$ for all $n<0, m\geq 0$.
Consequently, $P_{F_n,F_{n+m}}=r_1^{-m}z^m P_{F_n}$ and $P_{F_n,F_{n-m}}=r_1^{-m} (z^*)^m P_{F_n}$.
Hence, $P_{F_n, F_k} \in \mathcal I$ for all $n,k \leq -1$.
Moreover, $P_{F_n,E_k}=r_1^nz^mP_{F_n}$ and so $P_{F_n,E_k}\in\mathcal I$ for $k\geq 0$, $n\leq -1$.
In fact, it can be easily verif\/ied that, $P_{B_j, B_i}^*=P_{B_i, B_j}$.
This would mean that $P_{E_k,F_n}$ also belong to $\mathcal I$.

Similar calculations show that $P_{G_n, G_{n+m}}= r_2^{-m}P_{G_{n+m}}z^m P_{G_n}$ for $n+m<0$ and $P_{G_n, G_{n-m}}=
r_2^{-m}P_{G_{n-m}}(z^*)^m P_{G_n}$ for $n<0, m\geq 0$.
Collectively, these imply that $P_{G_n, G_k} \in \mathcal I$ for all \mbox{$n,k \leq -1$}.
Also $P_{G_n, E_{n+m}}= r_2^{m}P_{E_{n+m}}z^m P_{G_n}$ for $n+m \geq 0$.
Hence $P_{G_n, E_k}$ and $P_{E_k,G_n}$ belong to $\in \mathcal I$ for all $k \geq 0$, $n \leq -1$.

Finally, we notice that $P_{G_n, F_{-m}}= r_1^{-m}r_2^n P_{F_{-m}}(z^*)^m P_{E_0}z^n P_{G_n}$, which shows that the
operators $P_{G_n, F_k} \in \mathcal I$ for all $n,k \leq -1$.

Consider now f\/inite rank operators which are f\/inite linear combinations of the one-dimensio\-nal~$P_{B_i, B_j}$ for~$B_i$,~$B_j$ in the basis for~$H$.
It is a~simple exercise in functional analysis to show that all compact operators are norm limits of such f\/inite rank
operators.
\end{proof}

To describe the structure of $C^*(z)$ we need to understand the commutative quotient $C^*(z)/\mathcal{K}$.
For the case of quantum pair of pants that structure and the idea of proof is very similar to the case of quantum
annulus.

Below we will show that the $C^*$-algebra $C^*(z)$ contains some inf\/inite-dimensional projections.
Those are obtained from the spectrum of~$zz^*$.
While tedious, the computation of the spectrum of~$zz^*$ is fairly straightforward and it amounts to studying a~(multi
parameter) system of two step dif\/ference equations with constant coef\/f\/icients.
The results of the computations are presented in the next three theorems.

We start with the computation of the pure-point spectrum.

\begin{Theorem}
\label{discspecthm}
The operator $zz^*$ has three eigenvalues: $1$~with eigenspace $\operatorname{span}\{E_n\}_{n\geq 0}$, $r_1^2$~with eigenspace
$\operatorname{span}\{F_n\}_{n< 0}$, and the simple eigenvalue $\frac{r_1^2(a^2-r_2^2-r_1^2)}{a^2-r_1^2}$.
\end{Theorem}

\begin{proof}
We study $(zz^*-\lambda)\varphi=0$.
Using Lemma~\ref{z_z_star_coefficients} we get the following system of equations
\begin{alignat}{3}
& (1-\lambda)e_n=0 && \text{for}\quad n\ge1,&\nonumber
\\
 & \big(r_1^2+ r_2^2-\lambda\big)e_0+ar_2g_{-1}=0,&&&\nonumber
\\
&\big(r_1^2-\lambda\big)f_n=0 &&\text{for}\quad n\le-1,& \label{kernelequ}
\\
& ar_2g_{n+1}+\big(a^2+r_2^2-\lambda\big)g_n+ar_2g_{n-1}=0 \qquad && \text{for}\quad n\le-2,&\nonumber
\\
 & ar_2e_0+\big(a^2+r_2^2-\lambda\big)g_{-1}+ar_2g_{-2}=0.&&&\nonumber
\end{alignat}

The f\/irst and the third equations in this system yield the eigenvalues~1,~$r_1^2$ and the eigenspaces
$\operatorname{span} \{E_n\}_{n\geq 0}$, $\operatorname{span}\{F_n\}_{n< 0}$ respectively.
For the fourth equation, which is a~two step linear recurrence with constant coef\/f\/icients, the characteristic equation~is
\begin{gather*}
ar_2x^2+\big(r_2^2+a^2-\lambda\big)x+ar_2= 0.
\end{gather*}
The discriminant~$\Delta$ of this equation is
\begin{gather*}
\Delta=\big(r_2^2+a^2-\lambda\big)^2-4a^2r_2^2=\left((a-r_2)^2-\lambda\right)\left((a+r_2)^2-\lambda\right),
\end{gather*}
and the roots are
\begin{gather*}
x_{\pm}=\frac{\lambda-r_2^2-a^2 \pm \sqrt{\Delta}}{2ar_2}.
\end{gather*}
Thus, the formal solution to the homogeneous equation is
\begin{gather}
\label{ghomsol}
g_n=c_+x_+^{n+1}+c_-x_-^{n+1},
\end{gather}
where $c_+$ and $c_-$ are arbitrary constants.
Notice that if we adopt the convention $g_0:=e_0$ in the last equation of~\eqref{kernelequ}, then this formula holds for
$n\le 0$.

There are three separate cases to consider: $0<\lambda<(r_2-a)^2$, $(r_2-a)^2\leq \lambda\leq (r_2+a)^2$ and
$(r_2+a)^2<\lambda<1$.
If $0<\lambda<(r_2-a)^2$ we notice the following facts about the solutions.
Since $\lambda<(r_2-a)^2=r_2^2+a^2-2ar_2$ and both~$a$ and $r_2$ are positive, we have $\lambda-r_2^2-a^2<0$, and
since $\Delta>0$, we have that $x_-<0$.
Also notice that $x_+x_-=1$ from which it follows that $x_+<0$.
Since $x_-<x_+$ we have $x_-<-1$ and $-1<x_+<0$.
Thus, $|x_-|>1$ and $|x_+|<1$.

Since we need $g_n\in\ell^2({\mathbb Z}_{<0})$, and since $|x_+|<1$, we must have $c_+=0$.
In particular this implies that $g_0=x_-g_{-1}$.
The second equation of the system~\eqref{kernelequ} then becomes
\begin{gather*}
\big((r_1^2+r_2^2-\lambda)x_-+ar_2\big)g_{-1} =0.
\end{gather*}
If $(r_1^2+ r_2^2-\lambda)x_-+ar_2=0$, then using the relation $x_+x_- =1$ we obtain
\begin{gather*}
x_+=\frac{\lambda-r_1^2-r_2^2}{ar_2}.
\end{gather*}
On the other hand we also have
\begin{gather*}
x_+=\frac{\lambda-a^2-r_2^2+\sqrt{\Delta}}{ar_2}.
\end{gather*}
Setting these equal to each other and solving for~$\Delta$ we get $\Delta=(\lambda-2r_1^2-r_2^2+a^2)^2$.
Solving for~$\lambda$ yields
\begin{gather*}
\lambda=r_1^2\left(1-\frac{r_1^2}{a^2- r_1^2}\right)=\frac{r_1^2(a^2-r_2^2-r_1^2)}{a^2-r_1^2}.
\end{gather*}

First notice that due to the conditions on~$a$, $r_1$ and $r_2$ we have that $\lambda>0$ and $\lambda <r_1^2$ as
$r_1^2/(a^2-r_1^2)>0$.
Therefore this~$\lambda$ is in the interval $(0,(r_2-a)^2)$, and $c_-=g_{-1}$ is arbitrary.
Consequently $\lambda =\frac{r_1^2(a^2-r_2^2-r_1^2)}{a^2-r_1^2}$ is a~simple eigenvalue with an eigenvector
\begin{gather*}
\varphi_\lambda=\sum\limits_{n=-\infty}^0x_-^{n+1}G_n.
\end{gather*}

In the case $(r_2+a)^2<\lambda<1$, we have $\lambda-r_2^2-a^2=\lambda-(r_2+a)^2+2r_2a>0$, $\Delta>0$ and hence $x_+>0$.
Again from $x_+x_-=1$ we see that $x_+>1$ and $0<x_-<1$.
Consequently, in equation~\eqref{ghomsol} we must have $c_-=0$, which then implies that $g_0=x_+g_{-1}$.
The second equation of~\eqref{kernelequ} then becomes
\begin{gather*}
\big(\big(r_1^2+r_2^2-\lambda\big)x_++ar_2\big)g_{-1} =0.
\end{gather*}

Suppose there is~$\lambda$ such that $x_+(r_1^2+ r_2^2-\lambda)+ar_2=0$, which would then imply that
\begin{gather*}
x_+=\frac{ar_2}{\lambda-r_1^2-r_2^2}
\end{gather*}
and, since $x_+>1$, we must have
\begin{gather*}
\lambda < r_1^2+ar_2+r_2^2 < a^2+2ar_2+r_2^2=(r_2+a)^2.
\end{gather*}
This is a~contradiction.
Consequently $g_{n}=0$ for every~$n$ and there are no eigenvectors in this case.

When $(r_2-a)^2\leq\lambda\leq (r_2+a)^2$ the discriminant $\Delta\leq 0$ and the two solutions $x_+$ and $x_-$ are
complex numbers conjugate to each other with absolute value equal to one.
Since we need $g_n\in\ell^2({\mathbb Z}_{<0})$, equation~\eqref{ghomsol} implies that $g_{n}=0$ and there are again no
eigenvectors in this case.
\end{proof}

The second, and the most technical step in the calculation of the spectrum of $zz^*$ is the calculation of the inverse
of $zz^*-\lambda$.
Norm estimates of the inverse provide insight about the resolvent set of $zz^*$.

\begin{Lemma}
\label{resolvent_z_z_star}
The operator $(zz^*-\lambda)^{-1}$ is bounded for~$\lambda$ not an eigenvalue and
$\lambda\in(0,(r_2-a)^2)\cup((r_2+a)^2,1)$.
\end{Lemma}

\begin{proof}
We study $(zz^*-\lambda)\varphi=\tilde{\varphi}$ using coordinates
\begin{gather*}
\varphi=\sum\limits_{n=0}^\infty e_nE_n+\sum\limits_{n=-\infty}^{-1}(f_nF_n+g_nG_n)
\qquad
\text{and}
\qquad
\tilde{\varphi}=\sum\limits_{n=0}^\infty \tilde{e}_nE_n+\sum\limits_{n=-\infty}^{-1}\big(\tilde{f}_nF_n+\tilde{g}_nG_n\big).
\end{gather*}

Using Lemma~\ref{z_z_star_coefficients} we get the following system of equations
\begin{alignat}{3}
&  (1-\lambda)e_n=\tilde{e}_n && \text{for}\quad n\ge1,&\nonumber
\\
&\big(r_1^2+ r_2^2-\lambda\big)g_0+ar_2g_{-1}=\tilde{g}_0,&&&\nonumber
\\
&\big(r_1^2-\lambda\big)f_n=\tilde{f}_n && \text{for}\quad n\le-1,&\nonumber
\\
& ar_2g_{n+1}+\big(a^2+r_2^2-\lambda\big)g_n+ar_2g_{n-1}=\tilde{g}_n \qquad & &\text{for}\quad n\le-1.& \label{system_resolvent_z_z_star}
\end{alignat}
For the purpose of this proof we have introduced the notation $g_0:=e_0$ and $\tilde{g}_0:=\tilde{e}_0$.

The f\/irst and the third equations in system~\eqref{system_resolvent_z_z_star} can be solved directly
\begin{alignat*}{3}
& e_n=\frac{1}{1-\lambda}\tilde{e}_n \qquad && \text{for}\quad n\ge1
\quad
\text{and}
\quad
\lambda\neq 1, &
\\
& f_n=\frac{1}{r_1^2-\lambda}\tilde{f}_n \qquad && \text{for}\quad n\le-1
\quad
\text{and}
\quad
\lambda\neq r_1^2. &
\end{alignat*}

As for the fourth equation in system~\eqref{system_resolvent_z_z_star}, we start with the case $0<\lambda<(r_2-a)^2$.
From the proof of Theorem~\ref{discspecthm} we have that $|x_-|>1$ and $|x_+|<1$,
where $x_{\pm}=\frac{\lambda-r_2^2- a^2 \pm \sqrt{\Delta}}{2ar_2}$
are the solution of the characteristic equation $ar_2x^2+(r_2^2+a^2-\lambda)x+ar_2= 0$
with the discriminant $\Delta=((a-r_2)^2-\lambda)((a+r_2)^2-\lambda)$.

We use the variation of parameters technique to solve this equation.
For the homogeneous component of the solution, we use~\eqref{ghomsol} to obtain $g_n=c_+x_+^{n+1}+c_-x_-^{n+1}$ with
$c_+$ and $c_-$ being constants to be determined.
We look for the the solutions of the non-homogeneous equation in the form: $g_n=A_nx_+^{n+1}+B_nx_-^{n+1}$.
Then the standard trick is to assume the f\/irst equation below to obtain the following system
\begin{alignat*}{3}
&  (A_n -A_{n-1})x_+^{n+1}+(B_n-B_{n-1})x_-^{n+1}=0 \qquad &&\text{for}\quad  n\leq 0,&
\\
& (A_n -A_{n-1})x_+^n+(B_n-B_{n-1})x_-^n=\frac{-\tilde{g}_n}{ar_2} \qquad && \text{for}\quad  n\leq -1. &
\end{alignat*}
In particular,
\begin{gather}
\label{acompini}
(A_0 -A_{-1})x_++(B_0-B_{-1})x_-=0.
\end{gather}
The solution of the above system is
\begin{gather}
\label{acompgen}
A_n-A_{n-1}=\frac{\tilde{g}_nx_-^{n+1}}{\sqrt{\Delta}},
\qquad
B_n-B_{n-1}=\frac{-\tilde{g}_nx_+^{n+1}}{\sqrt{\Delta}}.
\end{gather}
In solving these dif\/ference equations we pay attention to square summability, making sure we only consider convergent
expressions in powers of $x_\pm$.
This leads to the following special solution of the non-homogeneous equation
\begin{gather*}
A_n=\frac{1}{\sqrt{\Delta}}\sum\limits_{j=-\infty}^n\tilde{g}_jx_-^{j+1}
\qquad
\text{and}
\qquad
B_n=\frac{1}{\sqrt{\Delta}}\sum\limits_{j=n+1}^{-1}\tilde{g}_jx_+^{j+1}
\end{gather*}
with $B_{-1}=0$, $n\leq -1$.
Consequently the general solution is
\begin{gather*}
g_n=c_-x_-^{n+1}+\frac{1}{\sqrt{\Delta}}\left(\sum\limits_{j=n+1}^{-1}\tilde{g}_jx_+^{j-n}+\sum\limits_{j=-\infty}^n\tilde{g}_jx_-^{j-n}\right)
\end{gather*}
for $n<-1$, since we want $g_n\in\ell^2({\mathbb Z}_{<0})$, and since $|x_+|<1$ we must have $c_+=0$.
For $n=-1$ we get
\begin{gather*}
g_{-1}=c_-+\frac{1}{\sqrt{\Delta}}\sum\limits_{j=-\infty}^{-1}\tilde{g}_jx_-^{j+1},
\end{gather*}
and for $n=0$, using~\eqref{acompini}, we obtain
\begin{gather*}
g_{0}=c_- x_-+A_0x_++B_0x_-= c_- x_-+A_{-1}x_++B_{-1}x_-=c_- x_-+\frac{x_+}{\sqrt{\Delta}}\sum\limits_{j=-\infty}^{-1}\tilde{g}_jx_-^{j+1}.
\end{gather*}

Next we study the second equation of system~\eqref{system_resolvent_z_z_star}.
Substituting the formulas for $g_{0}$ and $g_{-1}$ we compute
\begin{gather*}
c_-=\frac{\tilde{g}_0}{x_-(r_1^2+r_2^2-\lambda)+ar_2}-\left(\frac{{x_+(r_1^2+r_2^2-\lambda)+ar_2}}{x_-(r_1^2+r_2^2-\lambda)+ar_2}\right)
\frac{1}{\sqrt{\Delta}}\sum\limits_{j=-\infty}^{-1}\tilde{g}_jx_-^{j+1}.
\end{gather*}
The above formulas give the unique solution of the equation $(zz^*-\lambda)\varphi=\tilde{\varphi}$, and hence def\/ine
the inverse operator $(zz^*-\lambda)^{-1}$.
We need to verify that this operator is bounded.
So we f\/irst def\/ine the following operator
\begin{gather*}
Q\tilde{g}=\!\sum\limits_{n=-\infty}^0\!\! \left(\!\frac{\tilde{g}_0}{x_-(r_1^2+r_2^2-\lambda)+ar_2}-
\!\left(\frac{{x_+(r_1^2+r_2^2-\lambda)+ar_2}}{x_-(r_1^2+r_2^2-\lambda)+ar_2}\right)\!
\frac{1}{\sqrt{\Delta}}\! \sum\limits_{j=-\infty}^{-1}\!\! \tilde{g}_jx_-^{j+1}\!\right)\!x_-^{n+1}G_n.
\end{gather*}
Notice that, since $|x_-|>1$,~$Q$ is a~bounded operator taking $\ell^2({\mathbb Z}_{\leq 0})$ to itself.
Here we have used the notation $G_0:= E_0$.
Additionally we will need the following four operators, written in components:
\begin{gather*}
\left(L_1g\right)_n =\frac{1}{\sqrt{\Delta}}\sum\limits_{j=-\infty}^{-1}x_-^{j-n}\chi\left(\frac{j}{n}\right)g_j: \
\ell^2({\mathbb Z}_{<0})\to \ell^2({\mathbb Z}_{<0}),
\\
\left(L_2g\right)_n =\frac{1}{\sqrt{\Delta}}\sum\limits_{j=-\infty}^{-1}x_+^{j-n}\chi\left(\frac{n+1}{j}\right)g_j: \
\ell^2({\mathbb Z}_{<0})\to \ell^2({\mathbb Z}_{<0}),
\\
\left(L_3g\right)_n =\frac{1}{\sqrt{\Delta}}\sum\limits_{j=-\infty}^{-1}x_+^{j-n}\chi\left(\frac{j}{n}\right)g_j: \
\ell^2({\mathbb Z}_{<0})\to \ell^2({\mathbb Z}_{<0}),
\\
\left(L_4g\right)_n =\frac{1}{\sqrt{\Delta}}\sum\limits_{j=-\infty}^{-1}x_-^{j-n}\chi\left(\frac{n+1}{j}\right)g_j: \
\ell^2({\mathbb Z}_{<0})\to \ell^2({\mathbb Z}_{<0}).
\end{gather*}

Using these operators we write
\begin{gather*}
g_n=\left(Q\tilde{g}\right)_n+\left(L_1\tilde{g}\right)_n+\left(L_2\tilde{g}\right)_n \qquad\text{for}\quad  n\leq -1
\end{gather*}
and
\begin{gather*}
g_0 =\left(Q\tilde{g}\right)_0+x_+\left(L_1\tilde{g}\right)_{-1}.
\end{gather*}

Then we use Lemma~\ref{schuryounginq} (Schur--Young inequality) to estimate the norms of $L_1$ and $L_2$.
We have
\begin{gather*}
\|L_1\|^2
\le \frac{1}{\Delta}\left(\underset{n\le-1}{\sup}|x_-|^{-n}\sum\limits_{j=-\infty}^n|x_-|^j\right)
\left(\underset{j\le-1}{\sup}|x_-|^j\sum\limits_{n=j}^{-1}|x_-|^{-n}\right)
\\
\phantom{\|L_1\|^2}{}
 =\frac{1}{\Delta}\left(\frac{1}{1-|x_-|^{-1}}\right)\left(\underset{j\le-1}{\sup}\frac{1-|x_-|^{j}}{1-|x_-|^{-1}}\right)
 =\frac{1}{\Delta}\frac{1}{(1-|x_-|^{-1})^2}.
\end{gather*}
The computation of norm of $L_2$ is similar.
Therefore $(zz^*-\lambda)^{-1}$ is bounded for $0<\lambda<(r_2-a)^2$.

The other case is when $(r_2+a)^2<\lambda<1$.
From the proof of Theorem~\ref{discspecthm} we have that $x_+>1$ and $x_-<1$.
Again using variation of parameters, we see that the particular solution of system~\eqref{acompgen} is given by:
\begin{gather*}
B_n=\frac{-1}{\sqrt{\Delta}}\sum\limits_{j=-\infty}^n\tilde{g}_jx_+^{j+1}
\qquad
\text{and}
\qquad
A_n=\frac{-1}{\sqrt{\Delta}}\sum\limits_{j=n+1}^{-1}\tilde{g}_jx_-^{j+1}
\end{gather*}
with $A_{-1}=0$, $n\leq -1$.
Consequently the general solution is
\begin{gather*}
g_n=c_+x_+^{n+1} -\frac{1}{\sqrt{\Delta}}\left(\sum\limits_{j=n+1}^{-1}\tilde{g}_jx_-^{j-n}+\sum\limits_{j=-\infty}^n\tilde{g}_jx_+^{j-n}\right)
\qquad
\text{for}
\quad
n\leq-1.
\end{gather*}
Since we require $g_n\in\ell^2({\mathbb Z}_{<0})$, we must have $c_-=0$.
For $n=-1$ we obtain
\begin{gather*}
g_{-1}=c_+-\frac{1}{\sqrt{\Delta}}\sum\limits_{j=-\infty}^{-1}\tilde{g}_jx_+^{j+1},
\end{gather*}
and for $n=0$, using~\eqref{acompini}, we have
\begin{gather*}
g_{0}=c_+ x_++A_0x_++B_0x_-= c_+ x_++A_{-1}x_++B_{-1}x_-=c_+ x_+-\frac{x_-}{\sqrt{\Delta}}\sum\limits_{j=-\infty}^{-1}\tilde{g}_jx_+^{j+1}.
\end{gather*}

Substituting the formulas for $g_{0}$ and $g_{-1}$ into the second equation of system~\eqref{system_resolvent_z_z_star}
we com\-pu\-te~$c_+$
\begin{gather*}
c_+=\frac{\tilde{g}_0}{x_+(r_1^2+r_2^2-\lambda)+ar_2}+\left(\frac{{x_-(r_1^2+r_2^2-\lambda)+ar_2}}{x_+(r_1^2
+r_2^2-\lambda)+ar_2}\right) \frac{1}{\sqrt{\Delta}}\sum\limits_{j=-\infty}^{-1}\tilde{g}_jx_+^{j+1}.
\end{gather*}
The formulas above def\/ine the inverse operator $(zz^*-\lambda)^{-1}$.
To verify that this operator is bounded we def\/ine the following operator
\begin{gather*}
R\tilde{g}=\!\sum\limits_{n=-\infty}^0\!\!\left(\!\frac{\tilde{g}_0}{x_+(r_1^2+r_2^2-\lambda)+ar_2}-
\left(\frac{{x_-(r_1^2+r_2^2-\lambda)+ar_2}}{x_+(r_1^2+r_2^2-\lambda)+ar_2}\right)\!
\frac{1}{\sqrt{\Delta}}\sum\limits_{j=-\infty}^{-1}\!\!\tilde{g}_jx_+^{j+1}\!\right)\!x_+^{n+1}G_n.
\end{gather*}
It is easy to see that~$R$ is a~bounded operator taking $\ell^2({\mathbb Z}_{\leq 0})$ to itself.
Then we can write
\begin{gather*}
g_n=\left(R\tilde{g}\right)_n-\left(L_3\tilde{g}\right)_n-\left(L_4\tilde{g}\right)_n \qquad\text{for}\quad  n\leq -1
\end{gather*}
and
\begin{gather*}
g_0 =\left(R\tilde{g}\right)_0-x_-\left(L_3\tilde{g}\right)_{-1}.
\end{gather*}

We use Lemma~\ref{schuryounginq} to estimate the norms of $L_3$ and $L_4$, in a~manner similar to $L_1$.
We omit the repetitive details.
This shows that $(zz^*-\lambda)^{-1}$ is bounded for $(r_2+a)^2<\lambda<1$.
\end{proof}

In the theorem below, we will see that the interval $[(r_2-a)^2,(r_2+a)^2]$ is the continuous part of the spectrum of
$zz^*$, completing its full description.

\begin{Theorem}
The spectrum of $zz^*$ is
\begin{gather*}
\sigma(zz^*)= \left\{\frac{r_1^2(a^2-r_2^2-r_1^2)}{a^2-r_1^2}\right\}\cup\big\{r_1^2\big\}\cup \big[(r_2-a)^2,(r_2+a)^2\big]\cup\{1\}.
\end{gather*}
\end{Theorem}

\begin{proof}
Since $zz^*$ is a~positive operator with norm 1, its spectrum must be a~closed subset of the interval $[0,1]$.
In Theorem~\ref{discspecthm} we computed the pure point spectrum of $zz^*$ while Lemma~\ref{resolvent_z_z_star}
identif\/ied intervals belonging to the resolvent set of $zz^*$.
So it remains to analyze the interval $[(r_2-a)^2,(r_2+a)^2]$.
We will show that if $\lambda\in((r_2-a)^2,(r_2+a)^2)$ then $\operatorname{Ran}(zz^*-\lambda)$, the range of $(zz^*-\lambda)$, is not
all of~$H$.

In the notation of system~\eqref{system_resolvent_z_z_star} consider $(zz^*-\lambda)\varphi=\tilde{\varphi}$ with
$\tilde{g}_n=0$ for $n\leq -1$ but $\tilde{g}_0=\tilde{e}_0 \ne 0$.
This leads to the equation
\begin{gather*}
ar_2g_{n+1}+\big(a^2+r_2^2-\lambda\big)g_n+ar_2g_{n-1}=0 \qquad\text{for}\quad  n\leq -1,
\end{gather*}
with the general solution $g_n=c_+x_+^{n+1}+c_-x_-^{n+1}$.
Since for $(r_2-a)^2\leq\lambda\leq (r_2+a)^2$, the two numbers $x_+$ and $x_-$ are complex conjugates each with
magnitude one; which follows from the arguments in the last part of Theorem~\ref{discspecthm}; we must have $g_n =0$ for
$n\leq 0$ for~$\varphi$ to be in~$H$.
This however contradicts the second equation of~\eqref{system_resolvent_z_z_star} with $\tilde{g}_0=\tilde{e}_0 \ne 0$.
Consequently there is no $\varphi\in H$ satisfying $(zz^*-\lambda)\varphi=\tilde{\varphi}$ for such $\tilde{\varphi}$
and $\operatorname{Ran}(zz^*-\lambda)$ is not~$H$.
\end{proof}

Let $P_E$, $P_F$, and $P_G$ be the orthogonal projections onto the inf\/inite-dimensional span of $\{E_n\}_{n\geq 0}$,
$\{F_n\}_{n< 0}$ and $\{G_n\}_{n< 0}$ respectively.

\begin{Proposition}
The projections $P_E$, $P_F$, and $P_G$ belong to the noncommutative pair of pants, i.e., they are all in $C^*(z)$.
\end{Proposition}

\begin{proof}
Since Lemma~\ref{resolvent_z_z_star} implies that $r_1^2$ is an isolated eigenvalue of $zz^*$, there exists a~continuous
real valued function~$f$ so that $f(r_1^2)=1$ and~$f$ is zero on the rest of the spectrum of~$zz^*$.
By functional calculus we have $f(zz^*)=P_F$ and so $P_F\in C^*(z)$.
Similarly, since~$1$ is an isolated eigenvalue and we already know that $P_{E_0}\in C^*(z)$ we get $P_E\in C^*(z)$ as
well.
Since $P_G=I-P_E-P_F$, it then follows that $P_G\in C^*(z)$.
\end{proof}

We can decompose the Hilbert space~$H$ into $H \cong H_E \oplus H_F \oplus H_G$, where $H_E\cong\ell^2({\mathbb
Z}_{\ge0})$, $H_F\cong\ell^2({\mathbb Z}_{<0})$, and $H_G\cong\ell^2({\mathbb Z}_{<0})$ are the Hilbert spaces with
basis elements $E_n$, $F_n$, and $G_n$ respectively.
Since $\ell^2({\mathbb Z}_{\ge0})$ is a~subspace of $\ell^2({\mathbb Z})$, which can be identif\/ied with $L^2(S^1)$ via
the Fourier transform, we can view $H_E$ as a~subspace of $L^2(S^1)$.
More precisely, if $B_n=\zeta^n$, $n\in{\mathbb Z}$ is the standard basis in $L^2(S^1)$, then $H_E$ is identif\/ied with
the subspace $\operatorname{span}\{B_n\}_{n\geq 0}$ via
\begin{gather*}
E_n\mapsto B_n.
\end{gather*}

Def\/ine $P_{\ge0}: L^2(S^1)\to H_E$ to be the projection onto $H_E$.
For a~$\varphi\in C(S^1)$ we def\/ine $T_E(\varphi):H_E \to H_E$ by $T_E(\varphi)=P_{\ge0}M(\varphi)$ where
$M(\varphi):L^2(S^1)\to L^2(S^1)$ is the multiplication operator by~$\varphi$.
In particular for $\varphi(\zeta)=\zeta$ we have{\samepage
\begin{gather}
\label{teone}
T_E(\zeta)E_n=E_{n+1},
\end{gather}
the unilateral shift.}

Similarly we identify $H_F$ and $H_G$ with the subspace $\operatorname{span}\{B_n\}_{n< 0}$ in $L^2(S^1)$, and let $P_{<0}:L^2(S^1)\to
H_F$ and $P_{<0}':L^2(S^1)\to H_G$ be the orthogonal projections onto $H_F$ and $H_G$ respectively.
Then for a~$\varphi\in C(S^1)$, we def\/ine $T_F(\varphi):H_F \to H_F$ by $T_F(\varphi)=P_{<0}M(\varphi)$ and
$T_G(\varphi):H_G\to H_G$ by $T_G(\varphi)=P_{<0}'M(\varphi)$ respectively.
We have
\begin{gather}
\label{tfonf}
T_F(\zeta)F_n=F_{n+1},
\qquad
T_G(\zeta)G_n=G_{n+1}
\qquad
\text{for}
\quad
n<-1
\end{gather}
and
\begin{gather}
\label{tgong}
T_F(\zeta)F_{-1}=0,
\qquad
T_G(\zeta)G_{-1}=0.
\end{gather}

The operators $T_E(\varphi)$, $T_F(\varphi)$, and $T_G(\varphi)$ may be viewed as Toeplitz operators and they will be
needed in proving the following result.

\begin{Theorem}
\label{quotient_algebra}
The quotient $C^*(z)/\mathcal{K}$ is isomorphic to $C(S^1)\oplus C(S^1)\oplus C(S^1)$.
\end{Theorem}

\begin{proof}
Using the above notation def\/ine
\begin{gather*}
T: \  C\big(S^1\big)\oplus C\big(S^1\big)\oplus C\big(S^1\big) \to C^*(z)/\mathcal{K}
\end{gather*}
by
\begin{gather}
\label{pants_toeplitz_def}
T(\varphi_1,\varphi_2,\varphi_3)=T_E(\varphi_1)P_E+T_F(\varphi_2)P_F+T_G(\varphi_3)P_G+\mathcal{K}
\end{gather}
for continuous functions $\varphi_1$, $\varphi_2$, and $\varphi_3$ on the unit circle.
To see that~$T$ is well def\/ined we need to show that $T(\varphi_1,\varphi_2,\varphi_3)$ is in $C^*(z)$.
We showed that $P_E\in C^*(z)$, and notice that $T_E(\zeta)\in C^*(z)$ because $T_E(\zeta)=zP_E$.
But Toeplitz operators $T_E(\varphi)$ can be uniformly approximated by polynomials in $T_E(\zeta)$ and its adjoint, and
so $T_E(\varphi_1)P_E\in C^*(z)$.
Similar arguments work for $T_F(\varphi_2)P_F$ and $T_G(\varphi_3)P_G$.

We verify that~$T$ in~\eqref{pants_toeplitz_def} is a~isomorphism between the two algebras.
First notice that equation~\eqref{pants_toeplitz_def} implies that~$T$ is continuous and linear.
Next we show that the kernel of~$T$ is trivial.
Consider the equation $T(\varphi_1,\varphi_2,\varphi_3)=0$ implying that $T_E(\varphi_1)P_E+T_F(\varphi_2)P_F+T_G(\varphi_3)P_G$ is compact.
Since $P_E$, $P_F$ and $P_G$ are orthogonal, $T_E(\varphi_1)P_E$, $T_F(\varphi_2)P_F$ and $T_G(\varphi_3)P_G$ must be
compact, and consequently $T_E(\varphi_1):H_E \to H_E$, $T_F(\varphi_2):H_F \to H_F$, and $T_G(\varphi_3):H_G \to H_G$
are compact.
By the proof of Theorem~\ref{toeplitz_disk}, it follows that $\varphi_1=\varphi_2=\varphi_3=0$ and thus the kernel
of~$T$ is trivial.

Next we show that~$T$ is a~homomorphism of algebras.
Consider the dif\/ference:
\begin{gather}
T(\varphi_1,\varphi_2,\varphi_3)T(\psi_1,\psi_2,\psi_3)-T(\varphi_1\psi_1,\varphi_2\psi_2,\varphi_3\psi_3)
\nonumber
\\
\qquad{}
=\left(T_E(\varphi_1)T_E(\psi_1)-T_E(\varphi_1\psi_1)\right)P_E
\nonumber
\\
\qquad
\phantom{=}{}
+\left(T_F(\varphi_2)T_F(\psi_2)-T_F(\varphi_2\psi_2)\right)P_F+\left(T_G(\varphi_3)T_G(\psi_3)-T_G(\varphi_3\psi_3)\right)P_G.
\label{homo_T}
\end{gather}
Since $T_E$, $T_F$, and $T_G$ are Toeplitz operators, the proof of Theorem~\ref{toeplitz_disk} implies that all three
dif\/ferences on the right hand side of equation~\eqref{homo_T} are compact operators.
Thus~$T$ is a~homomorphism between the two algebras.

To show that the range of~$T$ is dense we consider the dif\/ference $T(\zeta,r_1\zeta, r_2\zeta+a)-z$.
Using formulas~\eqref{teone},~\eqref{tfonf}, and~\eqref{tgong} we get $T_E(\zeta)E_n=E_{n+1}$, for $n\geq 0$,
$T_F(r_1\zeta)F_n=r_1F_{n+1}$, for $n<-1$, and $T_G(r_2\zeta+a)G_n=r_2G_{n+1}+aG_n$, for $n<-1$.
Observe that $T(\zeta,r_1\zeta, r_2\zeta+a)-z$ is not zero on $F_{-1}$ and $G_{-1}$ only and hence it is a~compact operator.
Thus we have constructed functions $\varphi_1$, $\varphi_2$, and $\varphi_3$ such that $T(\varphi_1,\varphi_2,\varphi_3)
= z$ in $C^*(z)/\mathcal{K}$.
Since the $C^*$-algebra $C^*(z)/\mathcal{K}$ is generated by (the class of)~$z$, the range of~$T$ is dense and since the
range of a~$C^*$-morphism must be closed,~$T$ is an isomorphism of algebras.
This completes the proof.
\end{proof}

Note that from Theorem~\ref{quotient_algebra} we get a~short exact sequence
\begin{gather*}
0\rightarrow \mathcal{K}\rightarrow C^*(z)\rightarrow C\big(S^1\big)\oplus C\big(S^1\big)\oplus C\big(S^1\big)\rightarrow 0.
\end{gather*}
We can compare this to the short exact sequence for the classical pair of pants
\begin{gather*}
0\rightarrow C_0(PP_{(a, r_1,r_2)})\rightarrow C(PP_{(a, r_1,r_2)})\rightarrow C\big(S^1\big)\oplus C\big(S^1\big)\oplus C\big(S^1\big)\rightarrow 0,
\end{gather*}
where $C_0(PP_{(a, r_1,r_2)})$ are the continuous functions on the pair of pants that vanish on the boundary.

\pdfbookmark[1]{References}{ref}
\LastPageEnding


\begin{thebibliography}{99}
\footnotesize \itemsep=0pt

\bibitem{Abr}
Abrahamse M.B., Toeplitz operators in multiply connected regions,
  \href{http://dx.doi.org/10.2307/2373633}{\textit{Amer.~J. Math.}} \textbf{96} (1974), 261--297.

\bibitem{AD}
Abrahamse M.B., Douglas R.G., Operators on multiply connected domains,
  \href{http://dx.doi.org/10.2307/20488743}{\textit{Proc. Roy. Irish Acad. Sect.~A}} \textbf{74} (1974), 135--141.

\bibitem{Ahlfors}
Ahlfors L.V., Complex analysis. An introduction to the theory of analytic
  functions of one complex variable, 3rd ed., \textit{International Series in Pure and
  Applied Mathematics}, McGraw-Hill Book Co., New York, 1978.

\bibitem{Cob}
Coburn L.A., Singular integral operators and {T}oeplitz operators on odd
  spheres, \href{http://dx.doi.org/10.1512/iumj.1974.23.23036}{\textit{Indiana Univ. Math.~J.}} \textbf{23} (1973), 433--439.

\bibitem{Conway}
Conway J.B., A course in operator theory, \textit{Graduate Studies in
  Mathematics}, Vol.~21, Amer. Math. Soc., Providence, RI, 2000.

\bibitem{HS}
Halmos P.R., Sunder V.S., Bounded integral operators on~{$L^{2}$} spaces,
  \textit{Ergebnisse der Mathematik und ihrer Grenzgebiete}, Vol.~96,
  Springer-Verlag, Berlin~-- New York, 1978.

\bibitem{KL}
Klimek S., Lesniewski A., Quantum {R}iemann surfaces. {I}.~{T}he unit disc,
  \href{http://dx.doi.org/10.1007/BF02099210}{\textit{Comm. Math. Phys.}} \textbf{146} (1992), 103--122.

\bibitem{KL1}
Klimek S., Lesniewski A., Quantum {R}iemann surfaces. {II}.~{T}he discrete
  series, \href{http://dx.doi.org/10.1007/BF00402676}{\textit{Lett. Math. Phys.}} \textbf{24} (1992), 125--139.

\bibitem{KL2}
Klimek S., Lesniewski A., A two-parameter quantum deformation of the unit disc,
  \href{http://dx.doi.org/10.1006/jfan.1993.1078}{\textit{J.~Funct. Anal.}} \textbf{115} (1993), 1--23.

\bibitem{KL3}
Klimek S., Lesniewski A., Quantum {R}iemann surfaces. {III}.~{T}he exceptional
  cases, \href{http://dx.doi.org/10.1007/BF00761123}{\textit{Lett. Math. Phys.}} \textbf{32} (1994), 45--61.

\bibitem{KL4}
Klimek S., Lesniewski A., Quantum {R}iemann surfaces for arbitrary {P}lanck's
  constant, \href{http://dx.doi.org/10.1063/1.531503}{\textit{J.~Math. Phys.}} \textbf{37} (1996), 2157--2165.

\bibitem{KM1}
Klimek S., McBride M., D-bar operators on quantum domains, \href{http://dx.doi.org/10.1007/s11040-010-9084-9}{\textit{Math. Phys.
  Anal. Geom.}} \textbf{13} (2010), 357--390, \href{http://arxiv.org/abs/1001.2216}{arXiv:1001.2216}.

\bibitem{KM2}
Klimek S., McBride M., A note on {D}irac operators on the quantum punctured
  disk, \href{http://dx.doi.org/10.3842/SIGMA.2010.056}{\textit{SIGMA}} \textbf{6} (2010), 056, 12~pages, \href{http://arxiv.org/abs/1003.5618}{arXiv:1003.5618}.

\bibitem{KM3}
Klimek S., McBride M., Classical limit of the d-bar operators on quantum
  domains, \href{http://dx.doi.org/10.1063/1.3633525}{\textit{J.~Math. Phys.}} \textbf{52} (2011), 093501, 16~pages,
  \href{http://arxiv.org/abs/1101.2645}{arXiv:1101.2645}.

\bibitem{KM4}
Klimek S., McBride M., A note on gluing {D}irac type operators on a mirror
  quantum two-sphere, \href{http://dx.doi.org/10.3842/SIGMA.2014.036}{\textit{SIGMA}} \textbf{10} (2014), 036, 15~pages,
  \href{http://arxiv.org/abs/1309.7096}{arXiv:1309.7096}.

\bibitem{M}
Markushevich A.I., Theory of functions of a~complex variable, Chelsea
  Publishing Co., New York, 2005.

\end{thebibliography}
\end{document}